# Singular trajectories in multi-input time-optimal problems: Application to controlled mechanical systems

M. Chyba
Dept. of Mathematics
379 Applied Sciences Building
University of Santa Cruz
CA 95064

N.E. Leonard[*]
Dept. of Mechanical and Aerospace Eng.
D-234 Engineering Quad.
Princeton University
Princeton, NJ 08544

E.D. Sontag[†]
Dept. of Mathematics
Hill Center
Rutgers University
Piscataway, NJ 08854-8019

October 28, 2018

**Abstract**

This paper addresses the time-optimal control problem for a class of control systems which includes controlled mechanical systems with possible dissipation terms. The Lie algebras associated with such mechanical systems enjoy certain special properties. These properties are explored and are used in conjunction with the Pontryagin maximum principle to determine the structure of singular extremals and, in particular, time-optimal trajectories. The theory is illustrated with an application to a time-optimal problem for a class of underwater vehicles.

[*]Research partially supported by the Office of Naval Research under grant N00014-98-1-0649 and by the National Science Foundation under grant BES-9502477.
[†]Supported in part by US Air Force Grant F49620-01-1-0063

# 1 Introduction

In studying and designing controllers for physical processes, optimality of trajectories with respect to a given criterion is often a central concern. One might desire, for instance, to minimize the amount of fuel used by an airplane, or the energy spent by a given robotic system in order to reach a desired configuration. Similarly, one may wish to minimize flight time or operating time, since lengthy control efforts may be costly. In this paper, we focus on the time-optimal problem for a class of systems that includes controlled mechanical systems, such as rigid manipulators or underwater vehicles, for which we generally have more than one control. Controlled mechanical systems are of interest because of their practical significance in many applications and because of the special structure associated with the vector fields that define their dynamics. It is this structure that we make use of in this paper to better understand optimality of trajectories.

We consider both fully actuated and underactuated control systems. A first study on fully actuated conservative controlled mechanical systems was carried out in [10, 11]. As in that work, we base our study on the Pontryagin maximum principle, which gives necessary conditions for trajectories of control systems to be optimal with respect to a criterion such as energy or time. The idea is that one may associate, to any given controlled mechanical system, a set of vector fields describing the control system. Then, from the special form of the Lie algebra generated by these vector fields, one may extract from the maximum principle information on the structure of optimal trajectories. More precisely, according to the maximum principle, a time-optimal trajectory can be lifted to the cotangent bundle of its phase space as a trajectory of a constrained Hamiltonian system. This trajectory, combined with the corresponding control, is called an extremal. When the pointwise constraints in the maximum principle are nontrivial, one has nonsingular trajectories. These lead to boundary-valued controls, determined by the signs of the associated switching functions. However, it is well-known that an optimal trajectory may well be singular; that is, switching functions may vanish identically along the trajectory. The characterization of such trajectories, which is the question addressed in this paper, is in general a highly nontrivial problem. See for instance [1], [8] and [14] for a systematic use of Lie-algebraic techniques to restrict the trajectories candidates to time optimality.

Our main theorem gives conditions for the existence of singular extremals and the existence of trajectories with an infinite number of switchings. For instance, we can deduce from that theorem, for controlled mechanical systems with external forces, the nonexistence of totally singular extremals in the fully actuated situation (this result was already proved in [11]). We also give results on the possibility of concatenations of singular trajectories. In the case of underwater vehicles, we use these results to show that certain trajectories are not optimal.

While we focus here on mechanical systems (with the possibility of dissipative forcing), we note that our main result can be applied to other types of systems including nonholonomic robots.

In this paper, we perform computations using coordinates. However, it is of interest to translate our conditions into a coordinate-free formulation and to further study the geometry of our results. A geometrical approach to studying optimal problems associated with mechanical systems has been initiated in [6].

This paper is organized as follows. In Section 2 we state the problem and focus on the special class of affine control systems represented by controlled mechanical systems with external forces. In Section 3, we first recall the necessary conditions from the maximum principle for a trajectory of an affine control system to be time-optimal and introduce the switching functions which are the key tool to characterize the time-optimal trajectories. Then, we state and prove the main result with a specific application to the case of controlled mechanical systems. Section 4 is dedicated to the study of an application to illustrate the use of our main theorem. We treat the time-optimal problem for underwater vehicles and describe the singular trajectories.

## 2 Statement of the problem

In mechanics and robotics, most systems are described by an affine control system on a smooth $n$-dimensional manifold $\mathcal{M}$ (for Lagrangian systems with configuration space $\mathcal{Q}$, phase space is $\mathcal{M} = T\mathcal{Q}$). An affine control system on $\mathcal{M}$ takes the form

$$\dot{x}(t) = f(x(t)) + \sum_{i=1}^{m} g_i(x(t))u_i(t), \quad x(t) \in \mathcal{M} \tag{1}$$

where $u : [0, T] \to \mathcal{U} \subset {I\!\!R}^m$ is a measurable bounded function called the *control*, $f$ is a smooth vector field on $\mathcal{M}$, and the smooth vector fields $g_i$ are assumed to be linearly independent. They are called the *control vector fields*. The vector field $f$ is called the *drift* and represents the dynamics of the system. For a given control $u(.)$ defined on a time interval $[0, T]$ and an initial condition $x(0) = x_0$, a solution $x(.)$ of (1) is an absolutely continuous function defined on a sub-interval of $[0, T]$.

In this paper we will make the following assumption on the control-value set:

**Assumption 1.** The domain of control is given by

$$\mathcal{U} = \{u \in {I\!\!R}^m; \ \alpha_i \leq u_i \leq \beta_i, \ i = 1, \cdots, m, \ \alpha_i < 0 < \beta_i\}. \tag{2}$$

Such constraints on the inputs are very natural as the control usually represents quantities such as acceleration, temperature, etc., that cannot take arbitrarily big values. In the sequel an *admissible control* will refer to a measurable bounded function $u(.)$ defined on some interval $[0, T]$ such that $u(t) \in \mathcal{U}$ for almost every $t$.

We are interested in physical processes that can be governed; hence, it is natural to ask for the best control with respect to a given performance criterion. In this paper we focus on the *time-optimal problem*. The initial state and the target state of the system will be assumed to be points in the phase space $\mathcal{M}$. To summarize, we shall consider the following problem: given $x_i, x_f \in \mathcal{M}$ as initial and final states, find an admissible control $u(.)$ such that the corresponding trajectory steers the system (1) from $x_i$ to $x_f$ in the *shortest time*.

A key tool in our study will be the Lie algebra generated by the vector fields describing the system.

**Definition 1** *Let $X, Y$ be two smooth vector fields in $\mathcal{M}$. Their **Lie bracket** denoted by $[X, Y]$ is given by $[X, Y] = XY - YX$. The set of vector fields in $\mathcal{M}$ endowed with the Lie bracket operation is a **Lie algebra** denoted by $V^\infty(\mathcal{M})$. In the sequel $\mathcal{L}_G^f$ will denote the smallest sub-algebra of $V^\infty(\mathcal{M})$ that contains $f, g_1, \cdots, g_m$.*

## 2.1 Controlled mechanical systems

Controlled mechanical systems arise as follows. One starts with the specification of a Lagrangian of the form "kinetic minus potential energy". The kinetic energy is given by a Riemannian metric on the configuration manifold $\mathcal{Q}$. If we denote by $q = (q_1, \cdots, q_r)$ the local coordinates of the configuration manifold ($r$ is the number of degrees of freedom of the system), the Lagrangian defined on the phase space, $L : T\mathcal{Q} \to \mathbb{R}$, has the following expression in the local coordinates $(q, \dot{q})$ on $T\mathcal{Q}$:

$$L(q, \dot{q}) = \frac{1}{2}\dot{q}^t M(q)\dot{q} - V(q),$$

where $M(q)$ is a symmetric positive definite $r \times r$ matrix and $V(q)$ is a scalar function, the potential energy. We assume $M$ and $V$ are smooth with respect to $q$. The unforced equations of motion are given by the Euler-Lagrange equations and can be written in the form

$$M(q)\ddot{q} + C(q, \dot{q})\dot{q} + \frac{\partial V}{\partial q}^t = 0$$

where

$$C(q, \dot{q})\dot{q} = \frac{\partial}{\partial q}\left(M(q)\dot{q}\right)\dot{q} - \frac{\partial}{\partial q}\left(\frac{1}{2}\dot{q}^t M(q)\dot{q}\right). \tag{3}$$

The term $C(q, \dot{q})\dot{q}$, which is quadratic in $\dot{q}$, accounts for centrifugal and Coriolis forces. We also consider dissipation terms as well as other forces such as aerodynamics forces from fixed surfaces; we refer to these additional forces as *external forces* and collect them into an $r$-dimensional vector $D(q, \dot{q})$. To model the action of external controls, we introduce a smooth $r \times m$ matrix $Q(q)$ of rank $m$, and write the equations of motion in local coordinates in the form $Q(q)u = M(q)\ddot{q} + N(q, \dot{q})$, where

$$N(q, \dot{q}) = C(q, \dot{q})\dot{q} + \frac{\partial V}{\partial q}^t(q) - D(q, \dot{q})$$

is a smooth $r$-dimensional vector and $u(.)$ an admissible control. (Assuming that $Q$ depends only on configuration variables, and not velocities, means that we cannot include in the scope of our study some systems of interest, such as those involving control inputs corresponding to controlled surfaces on a moving body.) We will call the control system *conservative* (when $u = 0$) if $D \equiv 0$; this implies that the energy $E = \frac{1}{2}\dot{q}^t M \dot{q} + V$ is conserved. The above discussion leads to the following definition.

**Definition 2** *A **controlled mechanical system** is a system whose equations of motion are given in local coordinates by*

$$Q(q(t))u(t) = M(q(t))\ddot{q}(t) + N(q(t), \dot{q}(t)), \qquad q(t) \in \mathcal{Q} \qquad (4)$$

*where $M(q)$ is a symmetric positive definite $r \times r$ matrix, $N(q, \dot{q})$ is an $r$-dimensional vector and $Q(q)$ is an $r \times m$ matrix of rank $m$. The function $u : [0, T] \to \mathcal{U}$ is an admissible control as defined previously. All objects are assumed to be smooth. External forces (other than control forces) are included in the vector $N$. We say that (4) represents a controlled mechanical system **with quadratic external forces** if $N(q, \dot{q})$ is quadratic as a function of $\dot{q}$. The system is called **conservative** if $N(q, \dot{q}) = C(q, \dot{q})\dot{q} + \frac{\partial V}{\partial q}^t(q)$ for $C(q, \dot{q})\dot{q}$ given by (3) and for some smooth $V$. When $m = r$ (same number of input forces as number degrees of freedom), the system is fully actuated. It is underactuated if $m < r$.*

From now on, we will assume that the configuration manifold $\mathcal{Q}$ is the Euclidean space $\mathbb{R}^r$. As most of our results are local in nature, in applications they will apply to trajectories staying in a fixed chart. In the sequel, we introduce a smooth $n \times m-$matrix $G$ whose columns are the vector fields $g_i$. Hence $G$ is of rank $m$ and the system (1) can be rewritten as:

$$\dot{x}(t) = f(x(t)) + G(x(t))u(t), \quad x(t) \in \mathbb{R}^n \qquad (5)$$

A controlled mechanical system can be written as a control system of the form (5). When working with applications, see for instance Section 4, it is sometimes convenient to consider the following change of variables on the phase space from configuration variables and velocities $x = (q, \dot{q})^t$:

- $x_1 = \psi(q) \in \mathbb{R}^r$, where $q \mapsto \psi(q)$ is a smooth diffeomorphism;

- $x_2 = P(q)\dot{q} \in \mathbb{R}^r$, where $P(q)$ is a smooth $n \times n$ invertible matrix.

Given such $\psi$ and $P$, the equations of motion (4) lead us to an affine control system (1) where $x = (x_1, x_2)$ and $n = 2r$. More precisely, the drift is given by

$$f(x) = \begin{pmatrix} F(x_1)x_2 \\ \hat{f}(x_1, x_2) \end{pmatrix}_{2r \times 1} \qquad (6)$$

$$F(x_1) = \frac{d\psi}{dq}P^{-1}, \quad \hat{f}(x_1, x_2) = \dot{P}P^{-1}x_2 - PM^{-1}N(\psi^{-1}(x_1), P^{-1}x_2)$$

and the control vector fields by

$$G(x) = \begin{pmatrix} 0 \\ \hat{G}(x_1) \end{pmatrix}_{2r \times m}, \quad \hat{G}(x_1) = PM^{-1}Q \qquad (7)$$

the argument of the matrices being $\psi^{-1}(x_1)$. Note that $F$ is an invertible matrix of rank $r$. If the mechanical system is conservative or with quadratic external forces, $\hat{f}$ is of order 2 as a function of $x_2$. We denote by $(\hat{g}_i(x_1))_{i=1,\cdots m}$ the columns of the $r \times m$ matrix $\hat{G}$ and note that $g_i = (0, \hat{g}_i)^t$. Note that the vectors $\hat{g}_i$ are linearly independent, because the matrix $Q$ was assumed to be of rank $m$, so for the fully actuated situation $\hat{G}$ is an invertible $r \times r$ matrix.

### 2.1.1 Lie algebra associated to controlled mechanical systems

We already discussed the fact that one may naturally associate a Lie algebra of vector fields $\mathcal{L}_G^f$ to any controlled system. We now turn to a discussion of several special properties of $\mathcal{L}_G^f$ which hold whenever the control system in question is derived from a controlled mechanical system. The first essential property in our study is the commutativity condition satisfied by the vector fields $g_i$. Define $I = \{1, \cdots, m\}$.

**Lemma 1** *A controlled mechanical system satisfies the following commutativity condition:*

$$[g_i, g_j] \equiv 0, \qquad i, j = 1, \cdots, m. \qquad (8)$$

*Moreover, the vector fields $\{g_i, [f, g_i]\}_{i \in I}$ are pointwise linearly independent. If the system is fully actuated, the set of vector fields $\{g_i, [f, g_i]\}_{i \in I}$ generates, as a module over smooth functions, the set of all smooth vector fields on the phase space $I\!R^{2r}$.*

*Proof.* The columns of $G(x)$ are functions of the first $r$ state variables only, and their first $r$ components are zero. As a consequence, the relation $[g_i, g_j] = 0$ is verified for all $i, j \in I$. Moreover, computing the Lie brackets of the drift and the control vector fields we have:

$$[f, g_i] = \begin{pmatrix} -F(x_1)\hat{g}_i(x_1) \\ h_i(x_1, x_2) \end{pmatrix}$$

where $h_i$ is given by:

$$h_i(x_1, x_2) = \frac{\partial \hat{g}_i}{\partial x_1}(x_1)F(x_1)x_2 - \frac{\partial \hat{f}}{\partial x_2}(x_1, x_2)\hat{g}_i(x_1). \qquad (9)$$

Since $\hat{G}$ is a matrix of rank $m$ and $F$ is invertible, the family of vector fields $\{g_i, [f, g_i]\}_{i \in I}$ is linearly independent at every $x$.

Next, we make some preliminary remarks regarding the structure of Lie algebras generated by vector fields which have a polynomial structure in some of the variables; later we apply these remarks to controlled mechanical systems for which the external forces are polynomial as a function of $\dot{q}$.

By convention, we say that the polynomial "0" has degree $-1$. Assume we are given two vector fields $X, Y$ on $\mathbb{R}^n$ with the following properties: if we denote $x = (x_1, x_2) \in \mathbb{R}^{n_1} \times \mathbb{R}^{n_2}$, with $n_1 + n_2 = n$, then:

- As a function of $x_2$, the first $n_1$ coordinates of $X$ are polynomials of degree at most $a$, and the last $n_2$ of degree at most $b$.

- As a function of $x_2$, the first $n_1$ coordinates of $Y$ are polynomials of degree at most $c$, and the last $n_2$ of degree at most $d$.

Then, as a function of $x_2$, the first $n_1$ coordinates of the Lie bracket $[X, Y]$ are polynomial of degree at most $\max(a+c, a+d-1, b+c-1)$ and the last $n_2$ of degree at most $\max(b+c, a+d, b+d-1)$ with the convention that if $a$ (resp. $b, c$ or $d$) is negative then any term in these expressions containing $a$ (resp. $b, c$ or $d$) is -1, corresponding to the polynomial 0. In the sequel, $V_{n_1}^{a,b}$ will denote the set of vector fields in $\mathbb{R}^n$ whose first $n_1$ coordinates are polynomial in $x_2$ of degree at most $a$ and the last $n_2$ coordinates ($n_2 = n - n_1$) are polynomial in $x_2$ of degree at most $b$.

Of course, by adding more polynomial structure on the vector fields, we can obtain more conditions on the structure of the Lie brackets, but for our purposes the situation considered here is sufficient.

We now apply the above observation about degrees to controlled mechanical systems. In this case, we have $n = 2r$ and $n_1 = n_2 = r$.

**Lemma 2** *Assume that the external forces acting on a controlled mechanical system are polynomial, which means $N(q, \dot{q})$ is a polynomial as a function of $\dot{q}$. Then $f \in V_r^{1,b}$ for some $b \geq 2$ (in particular, $b = 2$ if the system is conservative), and as $g_i \in V_r^{-1,0}$, we have:*

$$ad_f^s g_i \in V_r^{(s-1)b-s+1, sb-s}, \qquad [g_j, ad_f^s g_i] \in V_r^{(s-1)b-s, sb-s-1} \qquad (10)$$

*for $i, j = 1 \cdots, m$. If the system is conservative or the external forces are of degree at most 2 as functions of $x_2$, then*

$$f \in V_r^{1,2}, \qquad ad_f^s g_i \in V_r^{s-1,s}, \forall s \geq 0, \qquad [g_j, [f, g_i]] \in V_r^{-1,0} \qquad (11)$$

*for $i, j = 1 \cdots, m$.*

**Corollary 1** *If a controlled mechanical system is fully actuated: $m = r$, the Lie brackets of the form $[g_j, [f, g_i]]$ belong to the module over the smooth functions on $\mathbb{R}^{2r}$ generated by the vector fields $g_i$.*

# 3 Time-optimality

We are interested in time-optimal trajectories for affine control systems:

$$\dot{x}(t) = f(x(t)) + G(x(t))u(t) \tag{12}$$

where $x(t) \in \mathcal{M} = \mathbb{R}^n$ (local study). Our study is based on the maximum principle, which provides necessary conditions for a trajectory to be optimal. For a general statement of the maximum principle, see [7] and [13] for instance. In the next section, we state the necessary conditions for our specific problem.

## 3.1 Maximum principle

Let $x_i, x_f \in \mathbb{R}^n$ be fixed, and suppose given an admissible time-optimal control $u : [0, T] \to \mathcal{U}$ such that the corresponding time-optimal trajectory $x(.)$ solution of (12) is defined on $[0, T]$ and steers the system from $x_i$ to $x_f$.

Using the maximum principle, we conclude that there exists an absolutely continuous vector $\lambda : [0, T] \to \mathbb{R}^n$, $\lambda(t) \neq 0$ for all $t$, such that the following conditions hold almost everywhere:

$$\dot{x}_j(t) = \frac{\partial H}{\partial \lambda_j}(\lambda(t), x(t), u(t)), \quad \dot{\lambda}_j(t) = -\frac{\partial H}{\partial x_j}(\lambda(t), x(t), u(t)) \tag{13}$$

for $j = 1, \cdots, n$, where $H(\lambda, x, u) = \lambda^t f(x) + \sum_{i=1}^m \lambda^t g_i(x) u_i$ is the Hamiltonian function, and the maximum condition holds:

$$H(\lambda(t), x(t), u(t)) = \max_{v \in \mathcal{U}} H(\lambda(t), x(t), v). \tag{14}$$

Moreover, the Hamiltonian is constant along the solutions of (13) and must satisfy $H(\lambda(t), x(t), u(t)) = \lambda_0$, $\lambda_0 \geq 0$.

**Definition 3** *A triple $(x, \lambda, u)$ which satisfies the maximum principle, in the sense just stated, is called an **extremal**. Its projection on the state space $x(\cdot)$ is said to be a **geodesic**, and the vector function $\lambda(\cdot)$ an **adjoint vector**. When the constant $\lambda_0$ is zero, the extremal is said to be **abnormal**.*

The maximum principle does not imply existence of optimal controls; for existence theorems the reader should refer to [4] for instance. As $\lambda(\cdot)$ is the solution of a linear differential equation, the condition "$\lambda(t) \neq 0$ for all $t$" is equivalent to: there exists $t_0 \in [0, T]$ such that $\lambda(t_0) \neq 0$.

**Definition 4** *The **switching functions** $\phi_i(.)$, $i = 1, \cdots, m$, along an extremal $(x, \lambda, u) : [0, T] \to \mathbb{R}^n \times \mathbb{R}^n \backslash \{0\} \times \mathcal{U}$ are defined by:*

$$\phi_i : [0, T] \to \mathbb{R}, \qquad \phi_i(t) = \lambda^t(t) g_i(x(t)). \tag{15}$$

*They are absolutely continuous functions.*

Clearly these functions play a crucial role in the study of time-optimal trajectories. Indeed, for systems of the form (12) and under Assumption 1 of Section 2 on the domain of control, the maximum condition (14) is equivalent to:

$$u_i(t) = \alpha_i \text{ if } \phi_i(t) < 0 \quad \text{and} \quad u_i(t) = \beta_i \text{ if } \phi_i(t) > 0, \qquad (16)$$

for $i = 1, \cdots, m$. Hence, the zeroes of the switching functions have to be carefully analyzed.

**Definition 5** *If there exists a nontrivial interval $[t_1, t_2] \subset [0, T]$ such that $\phi_i(t)$ is identically zero, the corresponding extremal is called $\mathbf{u_i}$-singular on $[t_1, t_2]$. The component $u_i$ of the control is then called singular on $[t_1, t_2]$. An extremal is called **totally singular** on $[t_1, t_2]$ if it is $u_i$-singular on $[t_1, t_2]$ for each $i$. The maximum principle implies that if $\phi_i(t) \neq 0$ for almost all $t \in [0, T]$ the component $u_i$ of the control is **bang-bang**, which means that it takes its values in $\{\alpha_i, \beta_i\}$ for almost every $t \in [0, T]$. Assume that the component $u_i$ of the control is bang-bang; then $t_s \in [0, T]$ is called a **switching time** for $u_i$ if, for each interval of the form $]t_s - \varepsilon, t_s + \varepsilon[ \cap [t_1, t_2]$, $\varepsilon > 0$, there is no constant $c$ such that $u_i(t) = c$ for almost all $t \in [t_1, t_2]$. An extremal is called $\mathbf{u_i}$-regular if $u_i$ is bang-bang with at most a finite number of switchings. An extremal is said to be **regular** if it is $u_i$-regular for all $i$.*

**Derivatives of the switching functions** The study of the zeroes of the switching functions is a key part of our analysis. Let us start by computing the first derivative: $\dot\phi_i(t)$. It is an easy verification that if $X$ is a smooth vector field in $\mathbb{R}^n$ and $(x, \lambda, u)$ an extremal, then the derivative of the absolutely continuous function $t \mapsto \lambda^t(t) X(x(t))$ with respect to $t$ is given by

$$\lambda^t(t)[f, X](x(t)) + \sum_{j=1}^{m} \lambda^t(t)[g_j, X](x(t)) u_j(t).$$

Hence, since $\phi_i(t) = \lambda^t(t) g_i(x(t))$, we have, almost everywhere in $t$:

$$\dot\phi_i(t) = \lambda^t(t)[f, g_i](x(t)) + \sum_{j=1}^{m} \lambda^t(t)[g_j, g_i](x(t)) u_j(t). \qquad (17)$$

Due to the fact that for generic affine control systems the Lie brackets $[g_i, g_j]$ do not commute, i.e. $[g_i, g_j] \neq 0$, the appearance of the control in the expression for $\dot\phi_i(.)$ means that this is in general merely a measurable function, and one cannot compute further derivatives. However, when considering controlled mechanical systems, Lemmas 1 and 2 give some additional structure on the Lie algebra $\mathcal{L}_G^f$. In particular, as controls disappear from (17), the derivatives of the switching functions are themselves absolutely continuous, and one can take further derivatives.

**Lemma 3** *Assume that the control system is derived from a controlled mechanical system. Then, if we write $\lambda = (\lambda_1, \lambda_2) \in \mathbb{R}^r \times \mathbb{R}^r$, the switching functions* depend only on $x_1$ and $\lambda_2$: $\phi_i(t) = \lambda_2^t(t)\hat{g}_i(x_1(t))$, $i = 1, \cdots, r$ where the vector fields $\hat{g}_i$ denote the columns of $\hat{G}$. The switching functions are continuously differentiable, their second derivatives exist and are measurable bounded functions. We have:

$$\dot{\phi}_i(t) = \lambda^t(t)[f, g_i](x(t)) \tag{18}$$
$$= -\lambda_1^t(t)F(x_1(t))\hat{g}_i(x_1(t)) + \lambda_2^t(t)h_i(x_1(t), x_2(t))$$

$$\ddot{\phi}_i(t) = \lambda^t(t)ad_f^2 g_i(x(t)) + \sum_{j=1}^m \lambda^t(t)[g_j, [f, g_i]](x(t))u_j(t). \tag{19}$$

where the functions $h_i$ are given by (9). In particular, if the system is conservative or the external forces are of degree at most 2 as a function of $x_2$, then $\dot{\phi}_i$ is linear as a function of $x_2$ and we have

$$\lambda^t(t)[g_j, [f, g_i]](x(t)) = \lambda_2^t(t)[g_j, \widehat{[f, g_i]}](x_1(t))$$

with $[g_j, \widehat{[f, g_i]}]$ an $r \times 1$ vector and $[g_j, [f, g_i]](x) = \begin{pmatrix} 0 \\ [g_j, \widehat{[f, g_i]}](x_1) \end{pmatrix}$.

If the system is fully actuated, the second derivative can be written as

$$\ddot{\phi}_i(t) = \lambda^t(t)ad_f^2 g_i(x(t)) + \sum_{j=1}^r \Big(\sum_{k=1}^r \alpha_{ij}^k(x(t))\phi_k(t)\Big)u_j(t) \tag{20}$$

where the $\alpha_{ij}^k$ are smooth functions defined by the relations

$$[g_j, [f, g_i]](x) = \sum_{k=1}^r \alpha_{ij}^k(x)g_k(x_1).$$

When the external forces are of degree at most 2 as functions of $x_2$, the coefficients $\alpha_{ij}^k$ depend only on $x_1$.

*Proof.* For a controlled mechanical system with external forces, $G$ is given by (7): $G(x) = \begin{pmatrix} 0 \\ \hat{G}(x_1) \end{pmatrix}$ where $\hat{G}$ is a $r \times m$ matrix depending on the $r$ first variables $x_1$ only. It follows that if $\lambda = (\lambda_1, \lambda_2) \in \mathbb{R}^r \times \mathbb{R}^r$ denotes the adjoint vector, we have $\lambda^t(t)g_i(x(t)) = \lambda_2^t(t)\hat{g}_i(x_1(t))$. From Lemma 1, the Lie brackets of the vector fields $g_i$ vanish, and then (17) becomes $\dot{\phi}_i(t) = \lambda^t(t)[f, g_i](x(t))$ and is an absolutely continuous function. Computing its derivative we obtain (19) which is a measurable bounded function. Moreover, if $N(q, \dot{q})$ is of degree at most 2 as a function of $\dot{q}$, from Lemma 2 we have $[g_j, [f, g_i]] \in V_r^{-1,0}$ for all $i, j = 1, \cdots, m$. If the system is fully actuated, then there exist smooth functions $\alpha_{ij}^k$ such that $[g_j, [f, g_i]](x) = \sum_{i=1}^r \alpha_{ij}^k(x_1)g_k(x_1)$, see Lemma 1. The rest of the proof follows easily.

A complete study of the time-optimal problem includes a complete characterization of the singular trajectories, an understanding of what concatenations of regular and singular extremals are possible, determination of a bound on the number of switching times (if one exists), proof of optimality of a sub-family of extremals (an optimal control between two given states is not necessarily unique for nonlinear systems), discussion of the regularity of an optimal synthesis, and so forth. Each of these questions is highly nontrivial.

For processes described by linear control systems and under controllability and normality conditions, see [4] or [9] for instance, existence and uniqueness theorems for time-optimal controllers are well-known and classical. In fact, in that case the conditions of the maximum principle are necessary and sufficient, and all time-optimal trajectories are regular. This leads to an optimal synthesis where the optimal control is given as a function of the state $x \in I\!R^n$. When the system is nonlinear, as is the case for most controlled mechanical systems, the time-optimal problem becomes much harder. In this paper, we are mainly concerned with the existence and characterization of singular extremals. Indeed, it is well known that an optimal control may well contain singular pieces. To analyze this question, we make use of geometric (Lie-theoretic) techniques. The power of such tools in the analysis of nonlinear time-optimal problems has been already illustrated by many authors; a good example and source of references is the paper [14] dealing with a car-like example with a polyhedral control-value set.

In the next section, we state the main result of this paper and its application to mechanical systems.

## 3.2 Theorem

Our main theorem deals with the existence of singular extremals and accumulation points of zeroes for the switching functions. It gives sufficient conditions on the Lie brackets involving the vector fields describing our affine control system under which we can draw conclusions. The application of this theorem to investigate time-optimality or lack thereof for system trajectories is illustrated in Section 4 for the underwater vehicle example.

In this section, $I$ denotes the set of the control indices: $I = \{1, \cdots, m\}$.

**Theorem 1** *Let $K_1, K_2 \subseteq I$ be such that $K_1 \cap K_2 = \emptyset$.*

*Assume that $(x, \lambda, u)$ is an extremal defined on $[0, T]$. If along the extremal there exists a finite sequence of nonempty sets $J_s \subseteq J'_{s-1} \subseteq J_{s-1} \cdots \subseteq J'_1 \subseteq J_1 \subseteq J_0 = K_1 \cup K_2$ such that for every $l \geq 1$ we have*

1.     $\forall j \in J_l$,
$$[g_k, ad_f^{l-1} g_j] \in Span\{ad_f^w g_v;\ w = 0, \cdots, l-1,\ v \in J_w\}$$
    $\forall k \in I$;

2.     $\forall j \in J'_l$,
$$[g_k, ad_f^{l-1} g_j] \in Span\{ad_f^w g_v;\ w = 0, \cdots, l-1,\ v \in \tilde{J}_w\}$$
    $\forall k \in I$ where $\tilde{J}_w = J_w \backslash \{J_w \cap K_2\}$;

3.
$$Span\{ad_f^w g_v;\ w=0,\cdots,s,\ v\in J_w\} = {I\!\!R}^n$$

*then either*

*a) there exists $i \in K_1$ such that the extremal is not $u_i$-singular;*

*or*

*b) there is no common accumulation point of zeroes of $\phi_i$, $i \in K_2$.*

If the subset $K_2$ is empty then condition 2 is equivalent to condition 1 and we have $J'_w = J_w$.

As we will see in Section 3.2.1 the following proposition gives information on the regularity of the nonsingular controls.

**Proposition 1** *Let $K_1, K_2 \subseteq I$, $K_1 \cap K_2 = \emptyset$. Assume $(x, \lambda, u)$ is an extremal defined on $[0, T]$ with the following properties:*

1. *It is $u_i$-singular for all $i \in K_1$.*

2. *There is a common accumulation point of zeroes for all $\phi_i$, $i \in K_2$ with corresponding time $t_0$.*

*If conditions 1,2 of Theorem 1 are satisfied for a finite sequence of nonempty sets $J_i, J'_i$ (as defined in Theorem 1) and $K \subseteq I \backslash \{K_1, K_2\}$ is such that*

$$Span\{g_u, ad_f^w g_v; w=1,\cdots,s, u\in K_1\cup K_2\cup K, v\in J_w\} = {I\!\!R}^n \quad (21)$$

*at $t = t_0$, then the switching functions $\phi_k$, $k \in K$ have no common zero at $t_0$. If $K_2$ is empty, then the switching functions $\phi_k$, $k \in K$ have no common zero along the whole extremal.*

Detailed proofs of Theorem 1 and Proposition 1 are given in Section 3.3. First, we study the consequences on Theorem 1 and Proposition 1 for the special structure of the Lie algebra $\mathcal{L}_G^f$ associated to controlled mechanical systems.

### 3.2.1 Application of Theorem 1 and Proposition 1 to controlled mechanical systems

**Lemma 4** *For a controlled mechanical system, conditions 1 and 2 of Theorem 1 are automatically satisfied for $l = 1$ (regardless of the sets $K_1, K_2$). For $l = 2$ condition 1 (resp. condition 2) is equivalent to*

$$[g_k, [f, g_j]] \in Span\{g_v;\ v \in J_0\ (\text{resp. } \tilde{J}_0)\}. \quad (22)$$

*Moreover,*
$$Span\{ad_f^w g_v;\ w=0,1, v\in J_w\} \text{ has dimension } d, \quad (23)$$

*where $d$ is the cardinality of $J_0 \cup J_1$.*

*Proof.* Condition 1 (resp. 2) for $l = 1$ is $[g_k, g_j] \in Span\{g_v; v \in J_0\}$ (resp. $v \in \tilde{J}_0$). It is automatically satisfied by the commutativity condition: $[g_k, g_j] = 0$, see Lemma 1. Relation (22) derives from the fact that the first $r$-coordinates of the vector fields $g_i$ and $[g_k, [f, g_j]]$ are zero (see Lemma 1) and that the $r$-first coordinates of the vector fields $[f, g_i]$ are linearly independent (see Lemma 1). Indeed, in any linear combination of the vector fields $[g_k, [f, g_j]]$ in terms of the vector fields $g_\ell$ and $[f, g_i]$, the coefficients of the latter must all vanish. From Lemma 1 we also know that the set consisting of all the vector fields $\{g_i, [f, g_j]\}$ is linearly independent, so we conclude that the subset in (23) also is.

**Proposition 2** *Assume given a fully actuated controlled mechanical system. Then:*

1. *Along an extremal we cannot have the same accumulation point of zeroes for all switching functions $\phi_i$, $i \in I$.*

2. *Let $k$ be a fixed element of $I$. Suppose along the extremal there exists $t_0 \in ]0, T[$ corresponding to a common accumulation of zeroes for all $\phi_i$, $i \neq k$. If $\phi_k(t_0) = 0$, then $\dot{\phi}_k(t_0) \neq 0$ and if*

$$Span\{g_i, [f, g_j], ad_f^2 g_j;\ i, j \in I, j \neq k\} = I\!R^{2r} \qquad (24)$$

   *at $t = t_0$, the switching function $\phi_k$ cannot vanish at $t = t_0$.*

*Proof.* The first part of the proposition is proved by applying Theorem 1 to $K_1 = \emptyset$, $K_2 = I$. Indeed in that case, by taking the sequence $J_0 = J_1 = I$ conditions 1,2 and 3 are satisfied (see Lemma 4). As $K_1$ is empty we deduce that there is no accumulation point of zeroes for all switching functions $\phi_i$, $i \in K_2 = I$. Assume now that $k$ is a fixed element of $I$ and $(x, \lambda, u)$ an extremal with a common accumulation point of zeroes at $t = t_0$ for all $\phi_i$, $i \neq k$. Since the functions $\phi_i$ are continuously differentiable, this implies that

$$\phi_i(t_0) = \lambda^t(t_0) g_i(x(t_0)) = 0,\ \dot{\phi}_i(t_0) = \lambda^t(t_0)[f, g_i](x(t_0)) = 0 \qquad (25)$$

for all $i \in I$, $i \neq k$. If $\phi_k$ and its first derivative vanish at $t = t_0$, equations (25) are also satisfied for $i = k$. The vector fields $g_i, [f, g_i]$ are linearly independent; hence this yields a contradiction with the nonvanishing condition for the adjoint vector in the maximum principle.

If, insteaad, (24) is verified, then we can apply Proposition 1 with $K_1 = \emptyset, K_2 = I \backslash \{k\}$ and $K = \{k\}$, the sequence of finite sets being given by $J_0 = J_1 = J_1' = J_2 = I \backslash \{k\}$. We use the fact from Lemma 1 that the Lie brackets of the form $[g_k, [f, g_j]]$ are obtained as combination over $C^\infty(R^{2r})$ of the vector fields $g_i$.

An immediate consequence is the nonexistence of totally singular extremal for a fully actuated controlled mechanical system. This result was already proved for conservative controlled mechanical systems in [11].

**Proposition 3** *A fully actuated controlled mechanical system ($m = r$) does not have totally singular extremals. More precisely, if $k$ is a fixed element of $I = \{1, \cdots, r\}$ and $(x, \lambda, u)$ is $u_i$-singular for all $i \neq k$, the extremal is $u_k$-regular. Moreover, if one of these two following conditions is satisfied, then the component $u_k$ of the control is constant:*

**(i)** $Span\{g_i, [f, g_j], ad_f^2 g_j; \ i, j \in I, j \neq k\} = \mathbb{R}^{2r}$ *along the extremal.*

**(ii)** *The extremal is abnormal, and along this extremal we have*

$$Span\{f, g_i, [f, g_j]; \ i, j \in I, j \neq k\} = \mathbb{R}^{2r}.$$

*Computation of the singular controls: along a time interval such that $\phi_k$ does not vanish, the $r-1$-singular components of the control are computed from the following set of equations:*

$$\frac{\lambda^t ad_f^2 g_i(x) + \alpha_{ik}^k(x)\phi_k u_k}{\phi_k} = -\sum_{j \neq k} \alpha_{ij}^k(x) u_j, \qquad i \in I \backslash \{k\} \qquad (26)$$

*where the coefficients $\alpha_{ij}^k$ are given by (20) and under the hypothesis that $\det(\alpha_{ij}^k(x)) \neq 0, \ j \neq k$.*

*Proof.* The proposition is a corollary of Proposition 2. It can also be proved directly by using Theorem 1 with $K_1 = I\backslash\{k\}$, $K_2 = \{k\}$ and Proposition 1 with $K_1 = I\backslash\{k\}$, $K_2 = \emptyset$ and $K = \{k\}$. To prove that $u_k$ is constant if condition (ii) is satisfied we proceed as follows. Along an extremal that is $u_i$-singular for all $i \neq k$, the Hamiltonian function becomes

$$H(\lambda(t), x(t), u(t)) = \lambda^t(t) f(x(t)) + \phi_k(t) u_k(t).$$

If the extremal is abnormal, $H \equiv 0$. Thus, if there exists $t_0$ such that $\phi_k(t_0) = 0$ we have $\lambda^t(t_0) f(x(t_0)) = 0$ (indeed $u_k$ is essentially bounded which means $\phi_k(t) u_k(t) \to 0$ as $t \to t_0$ except at most along a set of measure zero and $\lambda^t f(x)$ is continuous). This and condition (ii) implies $\lambda(t_0) = 0$ which is a contradiction. For a fully actuated controlled mechanical system, the second derivative is given by (20). Using the fact that $\phi_i \equiv 0$ for $i \neq k$, we obtain formula (26).

We have already noticed, see Lemma 1, that for a fully actuated controlled mechanical system, the set of vector fields $\{g_i, [f, g_i]\}_{i \in I}$ generates the set of all smooth vector fields on the phase space $\mathbb{R}^{2r}$. Hence, conditions 1-3 of Theorem 1 can be expressed in term of Lie brackets of length $\leq 2$. This is illustrated in [3] on the underwater vehicle example. We also use in that paper Theorem 1 to study extremals with fewer than $r-1$ components of the control singular at the same time.

In case the system is underactuated ($m < r$), sufficient conditions for the nonexistence of totally singular extremals are given by Theorem 1 in the following simplified form.

**Proposition 4** *Let $k$ be a fixed element of $I = \{1, \cdots, m\}$. Assume the extremal $(x, \lambda, u)$ is $u_i$-singular for all $i \neq k$ on the interval $[0, T]$. If there exists a finite sequence of nonempty sets $J_s \subseteq \cdots \subseteq J_3 \subseteq J_2 \subseteq I$ and $J_1 = J_0 = I$ such that along the extremal*

$$[g_k, [f, g_j]] \in Span\{g_v;\ v \in J_0\}, \qquad \forall j \in J_2, k \in I$$

$$[g_k, ad_f^{l-1} g_j] \in Span\{ad_f^w g_v;\ w = 0, \cdots, l-1,\ v \in J_w\}, \forall k \in I, l = 3, \cdots, s$$

$$Span\{ad_f^w g_v;\ w = 0, \cdots, s,\ v \in J_w\} = I\!\!R^{2n},$$

*then the extremal cannot be totally singular.*

To obtain information on the nonsingular component of the control we refer to Proposition 1. Proposition 4 will be illustrated on an underwater vehicle example in Section 4.

## 3.3 Proofs of Theorem and Proposition 1

Let us first prove the following technical lemma.

**Lemma 5** *Let $h : [0, T] \to I\!\!R$ be an absolutely continuous function. Assume that its first derivative (defined almost everywhere) is of the form*

$$\dot{h}(t) = \alpha(t) + \beta(t),$$

*where $\alpha : [0, T] \to I\!\!R$ is a continuous function and $\beta(.)$ is essentially bounded. Then, if $h(.)$ has an accumulation point of zeroes at $t = t_0 \in\ ]0, T[$ such that $\beta(t) \to 0$ as $t \to t_0$ (except along a set of zero measure), we have $\alpha(t_0) = 0$.*

*Proof.* An absolutely continuous function $h(.)$ is differentiable at almost all points of its definition interval. Assume $\dot{h}(t) = \alpha(t) + \beta(t)$ with $\alpha, \beta$ as stated in the lemma, then $\dot{h} : [0, T] \to I\!\!R$ is a measurable bounded function. Suppose $h(.)$ has an accumulation point of zeroes at $t = t_0$ and $\beta(t) \to 0$ as $t$ tends to $t_0$. Assume $\alpha(t_0) = c > 0$. Then there exists $\varepsilon > 0$ such that $\frac{c}{2} < \alpha(t) + \beta(t)$ almost everywhere on the interval $]t_0 - \varepsilon, t_0 + \varepsilon[$. This is equivalent to saying that there exists a set $F \subset [0, T]$ of measure zero such that if $t \in ]t_0 - \varepsilon, t_0 + \varepsilon[ \backslash F$ then $\dot{h}(t) \geq \frac{c}{2}$. For $t > t_0$, $h(t) = h(t_0) + \int_{t_0}^t \dot{h}(\tau) d\tau$ and $h(t_0) = 0$. So, we have $h(t) \geq \frac{c}{2}(t - t_0) > 0$ (if $t < t_0$, we have $h(t) < 0$). Hence $h(.)$ is a strictly increasing function in a neighborhood of $t_0$ which contradicts the fact that $t_0$ is an accumulation point of zeroes. If we assume $c < 0$, then we obtain the same contradiction with $h(.)$ a strictly decreasing function.

Let us now prove Theorem 1.

### 3.3.1 Proof of Theorem 1

By assumption, there exists $K_1, K_2 \subseteq I$ and a sequence of nonempty sets $J_i, J_i'$ as stated in the theorem such that conditions 1-2 and 3 are satisfied. Assume the extremal $(x, \lambda, u)$ is $u_i$-singular for all $i \in K_1$ and has a common accumulation point of zeroes for all $\phi_i$, $i \in K_2$. We will prove that it leads to a contradiction with the necessary conditions of the maximum principle. In the rest of the proof, we will denote by $t_0 \in ]0,T[$ any fixed time corresponding to a common accumulation point of zeroes for $\phi_i$, $i \in K_2$. If $K_2$ is empty, $t_0$ just denotes an arbitrary time (in particular in that case formula (31) is true for all $t$)

Assume that the extremal is defined on the time interval $[0,T]$.

<u>CLAIM:</u> If $j \in J_l$, the corresponding switching function $\phi_j$ is of class $C^{l-1}$, the $l$-th derivative exists, it is measurable bounded, and we have the following relations:

- <u>Case $l = 0$:</u>

$$\phi_j(t) = \lambda^t(t)g_j(x(t)) = 0, \ \forall t \in [0,T], j \in K_1, \tag{27}$$

$$\phi_j(t_0) = \lambda^t(t_0)g_j(x(t_0)) = 0, \ j \in K_2. \tag{28}$$

- <u>Case $l \geq 1$:</u>
  For $c = 0, \cdots, l-1$:

$$\phi_j^{(c)}(t) = \lambda^t(t)ad_f^c g_j(x(t)) = 0, \ \forall t \in [0,T], j \in J_c \cap K_1, \tag{29}$$

$$\phi_j^{(c)}(t_0) = \lambda^t(t_0)ad_f^c g_j(x(t_0)) = 0, \ j \in J_c \cap K_2, \tag{30}$$

and

$$\lambda^t(t_0)ad_f^l g_j(x(t_0)) = 0, \ \forall j \in J_l. \tag{31}$$

Moreover, the derivatives may be expressed as follows:

$$\phi_j^{(s)}(t) = \lambda^t(t)ad_f^s g_j(x(t)), \ \forall t \in [0,T], j \in J'_c, \tag{32}$$

$$\phi_j^{(l)}(t) = \lambda^t(t)ad_f^l g_j(x(t)) + \sum_{k=1}^m \lambda^t(t)[g_k, ad_f^{l-1}g_j]u_k(t), \ j \in J'_{l-1}, \tag{33}$$

for almost every $t$.

*Proof of CLAIM.*
First assume <u>$l = 0$</u>. In this case, we have $J_0 = K_1 \cup K_2$. If the extremal is $u_i$-singular for all $i \in K_1$, the corresponding switching functions are identically zero on $[0,T]$: $\phi_i \equiv 0$. In terms of vector fields, the following equations must be satisfied for all $i \in K_1$:

$$\phi_i(t) = \lambda^t(t)g_i(x(t)) = 0, \forall t \in [0,T]. \tag{34}$$

This is exactly (27). Moreover, relation (34) evaluated at $t = t_0$ is satisfied for all $i \in K_2$ which is equivalent to (28).

Assume now <u>$l = 1$</u>. Notice that in this case equations (29) and (30) are exactly (27) and (28), which we just proved. The switching functions being identically zero for all $i \in K_1$, we have in particular that their first derivatives, which are measurable bounded functions, vanish almost everywhere:

$$\dot{\phi}_i(t) = \lambda^t(t)[f, g_i](x(t)) + \sum_{k=1}^m \lambda^t(t)[g_k, g_i](x(t))u_k(t) = 0 \tag{35}$$

for a.e. $t \in [0, T]$, $\forall i \in K_1$. From condition 1 of the theorem, if $j \in J_1$ we have

$$[g_k, g_j] \in Span\{g_v;\ v \in J_0\}, \qquad \forall k \in I.$$

Hence, a consequence of (27) and (28) is that $\lambda^t(t_0)[g_k, g_j](x(t_0)) = 0$ for all $j \in J_1, k \in I$. Equation (35) implies then that $\lambda^t(t_0)[f, g_j](x(t_0)) = 0$ for all $j \in J_1 \cap K_1$. Indeed, as the components of the control are essentially bounded functions the terms $\lambda^t(t)[g_k, g_j](x(t))u_k(t)$ tend to 0 as $t$ tends to $t_0$ except at most along a set of measure zero. Now, since $\dot{\phi}_j \equiv 0$ and $\lambda^t(t)[f, g_j](x(t))$ is continuous, the latter function must vanish $t_0$. This gives equation (31) for $j \in J_1 \cap K_1$. To prove that this equation is also true for $j \in J_1 \cap K_2$ we proceed as follows. Formula (17) for the first derivative of the switching functions holds for any $i \in I$, not only under singularity assumptions. Hence if $j \in J_1 \cap K_2$, we have an accumulation point of zeroes at $t_0$ and $\dot{\phi}_j(t)$ is given by (17) where $\lambda^t(t_0)[g_k, g_j](x(t_0)) = 0$ for all $k \in I$. The situation is covered by Lemma 5 with $h(t) = \phi_j(t)$, $\alpha(t) = \lambda^t(t)[f, g_j](x(t))$ and $\beta(t) = \lambda^t(t)[g_k, g_j](x(t))$. Hence we can conclude that $\alpha(t_0) = 0$ for all $j \in J_1 \cap K_2$.

Let us next deal with the case $\underline{l = 2}$; the generalization to arbitrary $l$ will be then clear. By definition $J_2 \subset J_1'$. Let us then first describe condition 2 for $l = 1$. If $j \in J_1'$ and as $J_0 \setminus \{J_0 \cap K_2\} = K_1$ we have

$$[g_k, g_j] \in Span\{g_v;\ v \in K_1\}.$$

Hence equation (27) implies that $\lambda^t(t)[g_k, g_j](x(t))$ is identically zero along the extremal, for all $j \in J_1', k \in I$. This means that for $j \in J_1'$, the first derivative of $\phi_j$ is absolutely continuous and is given by $\dot{\phi}_j(t) = \lambda^t(t)[f, g_j](x(t))$, $\forall t \in [0, T]$. We then can compute the second derivative for these switching functions and obtain:

$$\ddot{\phi}_j(t) = \lambda^t(t) ad_f^2 g_j(x(t)) + \sum_{k=1}^{m} \lambda^t(t)[g_k, [f, g_j]](x(t))u_k(t),\ j \in J_1'. \qquad (36)$$

This is a measurable bounded function. Now, if $j \in J_2$ we have also from condition 1 that

$$[g_k, [f, g_j]] \in Span\{g_{v_0}, [f, g_{v_1}];\ v_0 \in J_0, v_1 \in J_1\}.$$

Hence, from equations (29)-(31) for $c = 0$, $l = 1$, we obtain the relation $\lambda^t(t_0)[g_k, [f, g_j]](x(t_0)) = 0$ for all $j \in J_2, k \in I$. Using this last relation, (36), and Lemma 5, in a similar way to what was done in the case $l = 1$, we complete the proof of the claim for $l = 2$.

Finally, we sketch the general case. Let $3 \leq l \leq s$ where $s$ is as stated in the theorem and assume the claim is satisfied for any $0 \leq c \leq l$. Then under the assumptions of the theorem we will prove that the claim is true for $l$. As the claim is true for $0 \leq c \leq l$, the switching functions $\phi_j$, $j \in J_{l-1}$ are of class $C^{(l-2)}$, and their $l-1$-th derivatives exist and are measurable bounded functions given by (33). From condition 1, if $j \in J_l$ we have

$$[g_k, ad_f^{l-1} g_j] \in Span\{ad_f^w g_v;\ w = 0 \cdots, l-1,\ v \in J_w\}.$$

Using equations (29)-(30) for $0 \leq c \leq l-1$, we obtain the following relation: $\lambda^t(t_0)[g_k, ad_f^{l-1}g_j](x(t_0)) = 0$, for all $j \in J_l, k \in I$. We complete the proof the same way as done for the cases $l = 1$ and $l = 2$, appealing to Lemma 5.

Having proved the claim, we are almost done with the proof of Theorem 1. Indeed condition 3 and relations (31) satisfied for $l = 0, \cdots, s$ imply a contradiction with the fact that the adjoint vector cannot vanish along an extremal. Hence, either there exists $i \in K_1$ such that the extremal is not $u_i$-singular, or there is no common accumulation point of zeroes of $\phi_i$, $i \in K_2$.

### 3.3.2 Proof of Proposition 1

Following the proof of Theorem 1 we have for $l \geq 0$:

$$\lambda^t(t_0)ad_f^l g_j(x(t_0)) = 0 \qquad \forall j \in J_l. \tag{37}$$

Moreover, if we have at $t = t_0$ a common zero for the switching functions $\phi_k$, $k \in K$ the following relation is true:

$$\lambda^t(t_0)g_k(x(t_0)) = 0 \qquad \forall k \in K. \tag{38}$$

Equations (37)-(38) and the condition given by formula (21) imply $\lambda(t_0) = 0$. This contradicts the maximum principle. In the case where $K_2 = \emptyset$, equation (37) is satisfied for all $t$. This completes the proof.

## 3.4 Concatenation of singular extremals

In order to understand optimal strategies, it is essential to determine which concatenations of pieces of extremals are allowed. This problem is in general very difficult. The following is a result on the concatenation of singular extremals along an optimal trajectory.

**Proposition 5** *Let $S_1, S_2$ be nonempty subsets of $I$. For $a = 1, 2$, let $(x^a, \lambda^a, u^a)$ be a $u_i$-singular extremal for all $i \in S_a$ and defined on a time interval $[0, T_a]$. Assume $x^1(T_1) = x^2(0)$. If for $a = 1, 2$, along $(x^a, \lambda^a, u^a)$ there exists a sequence of nonempty sets $J_s^a \subseteq J_{s-1}^a \subseteq \cdots \subseteq J_1^a \subseteq J_0^a = S_a$ such that for every $l \geq 1$ we have:*

1. $\forall j \in J_l^a$,

$$[g_k, ad_f^{l-1}g_j] \in Span\{ad_f^w g_v;\ w = 0, \cdots, l-1,\ v \in J_w^c\}$$

$\forall k \in I$.

2.
$$Span\{ad_f^w g_v(0);\ w = 0, \cdots, s,\ v \in J_w^1 \cup J_w^2\} = \mathbb{R}^n.$$

*Then, the concatenation $(x, \lambda, u)$ obtained by*

$$(x(t), \lambda(t), u(t)) = \begin{cases} (x^1(t), \lambda^1(t), u^1(t)) & t \in [0, T_1] \\ (x^2(T_1 - t), \lambda^2(T_1 - t), u^2(T_1 - t)) & t \in [T_1, T_1 + T_2] \end{cases}$$

*is not an extremal.*

*Proof.* First note that if $(x^2(t), \lambda^2(t), u^2(t))$ is an extremal defined on $[0, T_2]$ so is $(x_2(T_1-t), \lambda_2(T_1-t), u_2(T_1-t))$ defined on $[T_1, T_1+T_2]$. Hence, let us assume $(x^2, \lambda^2, u^2)$ is defined on $[T_1, T_1 + T_2]$. As the adjoint vector is defined modulo a multiplicative factor, we can assume $\lambda^2(T_1)$ to be any arbitrary value, and, in particular, we can assume $\lambda^2(T_1) = \lambda^1(T_1)$. If the concatenation $(x, \lambda, u)$ is an extremal then $\lambda : [0, T_1 + T_2] \to I\!\!R^n \backslash \{0\}$ is an absolutely continuous function satisfying equations (13). Under condition 1, using a similar proof as for Theorem 1 (with $K_2 = \emptyset$) we can show that for $j \in J_l^a$:

$$\lambda^t(t) ad_f^c g_j(x(t)) = 0, \qquad c = 0, \cdots, l \qquad (39)$$

for every $t \in [0, T_1]$ if $a = 1$ and every $t \in [T_1, T_1 + T_2]$ if $a = 2$. Hence, equation (39) and Condition 2 imply $\lambda(T_1) = 0$ which contradicts the maximum principle.

Applying this result to underwater vehicles, we will prove that some specific trajectories in the fully actuated situation are not optimal. Before doing so, we first give a result for fully actuated controlled mechanical systems.

**Proposition 6** *Consider a fully actuated controlled mechanical system. Then, an optimal trajectory cannot be the concatenation of a $u_i$-singular extremal for all $i \in S_1 \subset I$ and of a $u_i$-singular extremal for all $i \in S_2 \subset I$ if $S_1 \cup S_2 = I$.*

*Proof.* This is a direct consequence of Proposition 5 (take the sequences $J_0^a = J_1^a = S_a$).

## 4 Application to underwater vehicles

Our application concerns the motion planning problem for underwater vehicles. For such vehicles, it is not only worthwhile to minimize the energy spent to carry out a desired motion, but one is often also concerned with the continuous power consumption of all devices on board, such as sensors and computers. Because of this last quantity ("hotel load"), the amount of time used to travel between two desired configurations becomes a minimization criteria for our performance requirements. This section contains those results obtained from Theorem 1 when applied to such systems. Details of the computations, proofs, and more information about time-optimal trajectories, can be all found in [3, 2]. In [2], we also include a brief discussion on energy minimization, and its relation with some recent discoveries on marime mammals diving.

The dynamics of underwater vehicles can be described with the equations of motion as follows, see [5] for more details. The position and orientation of the underwater vehicle are identified with the group of rigid-body motions in $I\!\!R^3$: $SE(3) = \{(R, b); \ R \in SO(3), b \in I\!\!R^3\}$. If we define $\Omega = (\Omega_1, \Omega_2, \Omega_3)$ and $v = (v_1, v_2, v_3)$ to be the angular and translational velocity of the vehicle in body coordinates (see figure 1), then the kinematic equations are

$$\dot{R} = R\hat{\Omega}, \qquad \hat{\Omega} = \begin{pmatrix} 0 & -\Omega_3 & \Omega_2 \\ \Omega_3 & 0 & -\Omega_1 \\ -\Omega_2 & \Omega_1 & 0 \end{pmatrix} \qquad (40)$$

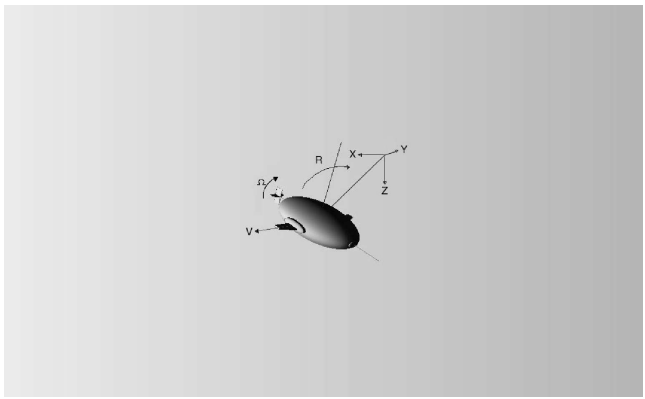

Figure 1: underwater vehicle

$$\dot{b} = Rv. \tag{41}$$

We begin by assuming the vehicle is submerged in an infinitely large volume of incompressible, irrotational and inviscid fluid at rest at infinity. Let us denote by $\Pi$ and $P$ the angular and linear components of the impulse (roughly equivalent to momentum) of the body-fluid system with respect to the body-frame. The Kirchhoff equations of motion for a rigid body in such an ideal fluid are given by

$$\dot{\Pi} = \Pi \times \Omega + P \times v + \tau, \;\; \dot{P} = P \times \Omega + \mathcal{F} \tag{42}$$

where $\tau$ and $\mathcal{F}$ are external torque and force vectors. $\tau$ and $\mathcal{F}$ can be used to include gravity, buoyancy and control forces as well as viscous forces such as lift and drag. Notice that $\Pi$ and $P$ can be computed from the total kinetic energy $T$ of the body fluid system: $T = \frac{1}{2}(\Omega^t J \Omega + 2\Omega^t D v + v^t M v)$, where $M$ and $J$ are respectively the body-fluid mass and inertia matrices. In [5], the Hamiltonian structure of the dynamics of the underwater vehicle is described.

We consider a neutrally buoyant, uniformly distributed, ellipsoidal vehicle moving in the vertical inertial plane and we neglect viscous effects. We denote by $(x, z)$ the absolute position of the vehicle where $x$ is the horizontal position and $z$ the vertical position. $\theta$ describes its orientation in this plane so that $q = (x, z, \theta)$, see figure 2. Given our assumptions, $M$ and $J$ are diagonal and $D = 0$ so that for our vehicle restricted to the plane $T = \frac{1}{2}(I\Omega^2 + m_1 v_1^2 + m_3 v_3^2)$ where $I$ is the body-fluid moment of inertia in the plane and $m_1, m_3$ are body-fluid mass terms in the body horizontal and vertical directions, respectively. We assume that $m_1 \neq m_3$, i.e., the planar vehicle is not a circle. We choose the state vector to be $w = (x, z, \theta, v_1, v_3, \Omega)$. Here $\Omega$ is the scalar angular rate in the plane. The equations of motion of this mechanical system are

$$\begin{aligned}\dot{x} &= \cos\theta v_1 + \sin\theta v_3 \\ \dot{z} &= \cos\theta v_3 - \sin\theta v_1\end{aligned}$$

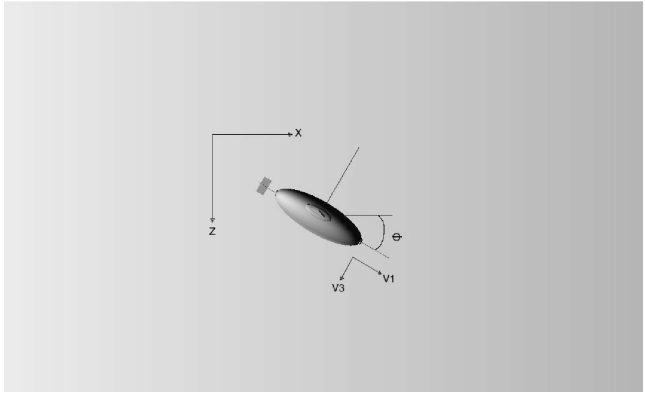

Figure 2: planar version

$$\begin{aligned}
\dot{\theta} &= \Omega \qquad (43)\\
\dot{v}_1 &= -v_3 \Omega \frac{m_3}{m_1}\\
\dot{v}_3 &= v_1 \Omega \frac{m_1}{m_3}\\
\dot{\Omega} &= v_1 v_3 \frac{m_3 - m_1}{I}
\end{aligned}$$

Let us consider the fully actuated case: the control vector is $u = (u_1, u_2, u_3)$ where $u_1$ is a force in the body 1-axis, $u_2$ is a force in the body 3-axis and $u_3$ is a pure torque in the plane. Accordingly, with the drift vector field $f$ given by the planar equations of motion (44), the input vector fields $g_i$ are given by

$$g_1 = \begin{pmatrix} 0 \\ 0 \\ 0 \\ \frac{1}{m_1} \\ 0 \\ 0 \end{pmatrix}, \quad g_2 = \begin{pmatrix} 0 \\ 0 \\ 0 \\ 0 \\ \frac{1}{m_3} \\ 0 \end{pmatrix}, \quad g_3 = \begin{pmatrix} 0 \\ 0 \\ 0 \\ 0 \\ 0 \\ \frac{1}{I} \end{pmatrix}.$$

With respect to the notations of Section 2.1 the change of variables on the phase space is given by: $\psi = Identity$ and $P^{-1}(w) = Q(w)$ where $Q(w)$ is a rotation matrix $\begin{pmatrix} R(\theta) & 0 \\ 0 & 1 \end{pmatrix}$, $R(\theta) = \begin{pmatrix} \cos\theta & \sin\theta \\ -\sin\theta & \cos\theta \end{pmatrix}$. This formulation of the problem (instead of considering the velocities expressed in the inertial frame coordinates) leads to simpler computations for the Lie brackets and to a nice geometric interpretation of the results.

The following result is an immediate consequence of Proposition 3.

**Proposition 7** *There is no totally singular extremal in the time-optimal problem for fully actuated underwater vehicles.*

As a particular consequence, if the underwater vehicle follows the equation of motion of the conservative mechanical system (which corresponds to the control identically zero), such a motion is not time-optimal. In fact, from Proposition 2 we have more information: along an extremal, there cannot exist any common accumulation point of zeroes for all switching functions.

Thus, we must first study extremals with at most two components of the control singular at the same time. From Proposition 3 we know that the nonsingular component of the control is regular and constant if certain Lie brackets conditions are satisfied. Let us describe more precisely what happens in the case of the fully actuated underwater vehicles. The Lie algebra generated by the Lie brackets of the vector fields $f, g_i$ is central to the characterizations of time optimal paths, and hence is the key to the proofs of the results which we are stating here. In order to avoid overloading with formulas, we refer the reader interested in by the computations and the structure of the Lie algebra to [3, 2]. Observe that, using notations and Lemmas of Section 2, we have the following facts: Our underwater vehicle is a conservative controlled mechanical system, hence there are no external forces and we have $f \in V_3^{1,2}$. It follows that $ad_f^s g_i \in V_3^{s-1,s}$ and $[g_j, ad_f^s g_i] \in V_3^{s-2,s-1}$ (this is verified in [3], where the Lie brackets are computed). The next proposition describes the 2-singular extremals.

**Proposition 8**

- Along a $u_1, u_2$-singular extremal, the component $u_3$ of the control is regular, and has at most one switching contained in the set $\{w;\ v_1 = v_3 = 0\}$. If there is one, then $v_1 \equiv v_3 \equiv 0$ along the extremal, and the singular components of the control are identically zero. The corresponding motion for the underwater vehicle is a pure rotation with angular velocity varying linearly.

- Let $i \in \{1, 2\}$ and $k$ be such that $k = 1$ if $i = 1$, and $k = 3$ if $i = 2$. Along a $u_i, u_3$-singular extremal, there is at most one switching for the nonsingular component of the control and this switching is contained in the set $\{w;\ v_k = \Omega = 0\}$. Moreover, if there is one switching, then $\Omega \equiv 0$ and $v_k \equiv 0$, and both singular components of the control are identically zero. The corresponding motions for the underwater vehicle are translations in the direction of the vertical (if $i = 1$) and horizontal (if $i = 2$) axis of the body frame coordinates, and the velocity of these translations depends linearly on $t$.

- Along an abnormal extremal with 2 singular components of the control, there is a switching at the time $t_s$ only if the underwater vehicle is at rest at $t = t_s$: $v_1(t_s) = v_3(t_s) = \Omega(t_s) = 0$.

In [3] we describe an algorithm to compute the singular components of the control when there is no switching and illustrate our results with simulations. The previous proposition shows that along a 2-singular extremal there is at most one switching for the nonsingular component of the control. This is an important result, indeed the existence of a uniform bound is not a consequence of the maximum principle. This is a well-known question in optimal control: in [12] for instance, the author provides an example of optimal trajectories with an infinite number of switchings. A crucial question is: is there a uniform bound on the number of switchings for all optimal trajectories? This is a highly nontrivial question, and we will try to give at least partial answers in forthcoming articles.

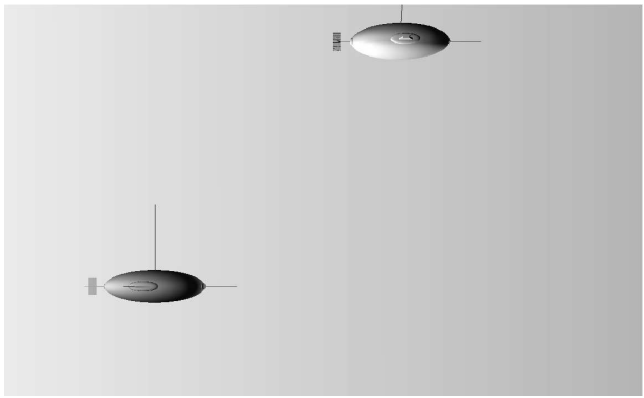

Figure 3: initial and final position

Conversely to Proposition 8, the following proposition states that the three basic motions: pure rotation,on and translations in the direction of the axis of the body frame coordinates, are 2-singular extremals.

**Proposition 9** *If along a trajectory we have $v_1 \equiv 0$ and $v_3 \equiv 0$, then the trajectory is a $u_1, u_2$-singular extremal with at most one switching. The corresponding motion is a pure rotation and the singular components of the control are identically zero: $u_1 \equiv u_2 \equiv 0$. If along a trajectory we have $\Omega \equiv 0$ and $v_1 \equiv 0$ (resp. $v_3 \equiv 0$), then the trajectory is a $u_1, u_3$-singular extremal (resp. $u_2, u_3$) with at most one switching for $u_2$ (resp. $u_1$). The corresponding motion is a translation in the direction of the vertical (resp. horizontal) axis of the body frame coordinates.*

We conjecture the time-optimality of these basic motions. THis would imply that one cannot restrict the study to bang-bang trajectories when dealing with time-optimality for such systems.

Let us use Proposition 6 to discuss the optimality of some specific trajectories. Consider the initial and final positions illustrated in figure 3. Assume the underwater vehicle to be at rest at these positions. On figure 4 and 5 we represent two different ways to connect these positions; these are concatenations of vertical and horizontal translations (the orientation is fixed along the motion). Using the result given in Proposition 6 on the concatenation of singular extremals in the fully actuated situation, we can prove that these trajectories are not optimal.

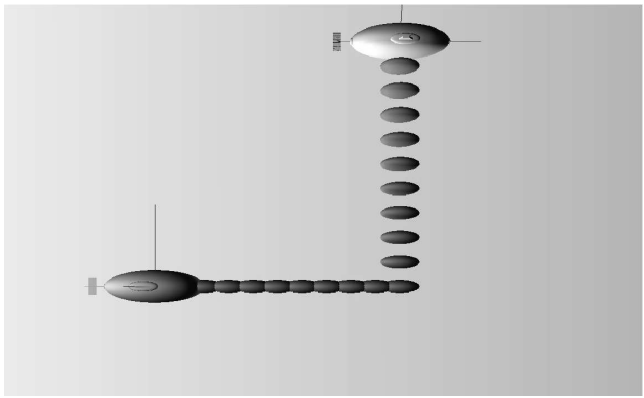

Figure 4: vertical + horizontal translation

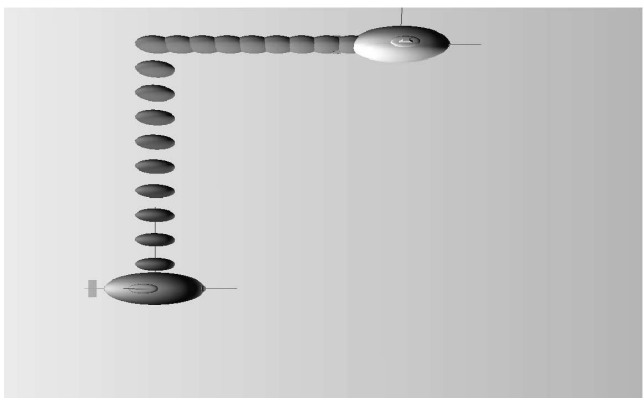

Figure 5: horizontal + vertical translation

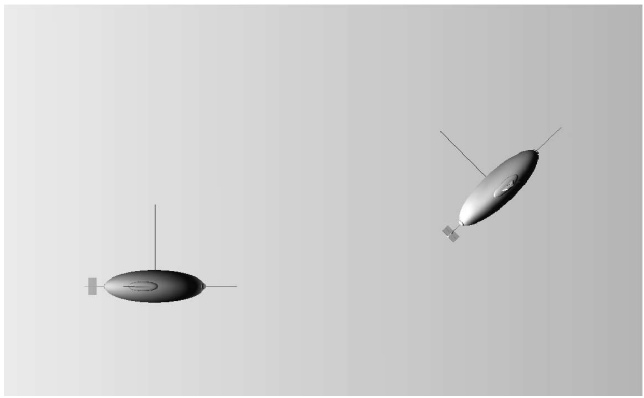

Figure 6: Initial and final configurations

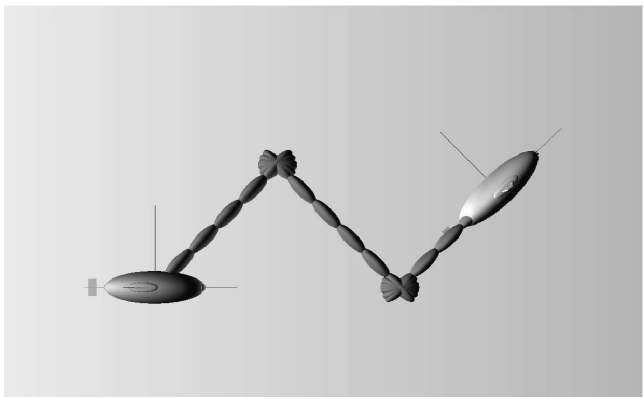

Figure 7: path 1

More generally, we have the following result. Let $w_i$ and $w_f$ arbitrary initial and final positions such that the velocities are all zero. We can reach these two positions by a motion formed by pure rotation and pure translations in the direction of the axis of the body frame coordinates. This is illustrated by figure 6 and 7. We know that these trajectories are not optimal. The nonoptimality of such trajectories is guaranteed by the fact that pure translational and rotational motions are necessarily 2-singular extremals, and by Proposition 6, which does not allow concatenations of 2-singular extremals in our case.

In order to obtain more information on the stucture of the time-optimal trajectories, we also need to study extremals with only one component of the control being singular, as well as trajectories with all components of the control bang-bang. This will be done in a forthcoming article using results presented in the previous sections. For instance, Proposition 2 gives information about extremals with one component singular only.

**Proposition 10** *Let $\{i, j, k\} \in \{1, 2, 3\}$ such that $i \neq j \neq k$. If the extremal is $u_i$-singular, then $\phi_j$ and $\phi_k$ cannot have a same accumulation point of zeroes. If there exists $t_0$ such that $\phi_j(t_0) = \phi_k(t_0) = 0$, we have $\dot{\phi}_j(t_0) \neq 0$ or $\dot{\phi}_k(t_0) \neq 0$.*

As we already mentionned in this paper, our main result does not apply exclusively to the fully actuated situation. In a forthcoming article, we will study underactuated situations. These situations will be motivated by a buoyancy driven underwater vehicle; see [2] for a first approach.